\documentclass[amsfonts]{amsart}
\usepackage[margin=1.25in]{geometry}

\def\UseBibLatex{1}

\makeatletter
\def\input@path{{styles/}}
\makeatother

\newcommand{\UsePackage}[1]{%
  \IfFileExists{styles/#1.sty}{%
      \usepackage{styles/#1}%
   }{%
      \IfFileExists{../styles/#1.sty}{%
         \usepackage{../styles/#1}%
      }{%
         \usepackage{#1}%
      }%
   }%
}

\providecommand{\BibLatexMode}[1]{}
\providecommand{\BibTexMode}[1]{}

\ifx\UseBibLatex\undefined%
  \renewcommand{\BibLatexMode}[1]{}
  \renewcommand{\BibTexMode}[1]{#1}
\else
  \renewcommand{\BibLatexMode}[1]{#1}
  \renewcommand{\BibTexMode}[1]{}
\fi

\BibLatexMode{%
   \usepackage[bibencoding=utf8,style=alphabetic,backend=biber]{biblatex}%
   \UsePackage{my_biblatex}%
}

\usepackage{amsthm}

\newtheorem*{thm*}{Theorem}

\newtheorem{thm}{Theorem}[section]
\newtheorem{cor}[thm]{Corollary}
\newtheorem{prop}[thm]{Proposition}
\newtheorem{lem}[thm]{Lemma}

\theoremstyle{definition}

\newtheorem{notn}[thm]{Notation}

\theoremstyle{remark}
\newtheorem{rem}[thm]{Remark}

\usepackage{amsmath,amssymb,tikz-cd,calligra,mathrsfs,biblatex, hyperref}
\newcommand{\Z}{\mathbb{Z}}

\DeclareMathOperator{\Imag}{Im}
\DeclareMathOperator{\Ker}{Ker}

\DeclareMathOperator{\rk}{rk}
\DeclareMathOperator{\Id}{Id}
\DeclareMathOperator{\Tr}{Tr}
\DeclareMathOperator{\Hom}{Hom}
\DeclareMathOperator{\SHom}{\mathop{\mathcal{H}\! \mathit{om}}}
\DeclareMathOperator{\End}{End}
\DeclareMathOperator{\Mat}{Mat}

\DeclareMathOperator{\GL}{GL}
\DeclareMathOperator{\Gr}{Gr}

\DeclareMathOperator{\Sym}{Sym}
\DeclareMathOperator{\Ind}{Ind}
\DeclareMathOperator{\Res}{Res}
\DeclareMathOperator{\Spec}{Spec}

\DeclareMathOperator{\pr}{pr}

\DeclareMathOperator{\Tot}{Tot}
\DeclareMathOperator{\gra}{graph}

\newcommand{\mcB}{\mathcal{B}}

\newcommand{\mcO}{\mathcal{O}}
\newcommand{\mcV}{\mathcal{V}}
\newcommand{\mcF}{\mathcal{F}}
\newcommand{\mcM}{\mathcal{M}}
\newcommand{\mcN}{\mathcal{N}}
\newcommand{\mcQ}{\mathcal{Q}}
\newcommand{\mcS}{\mathcal{S}}

\numberwithin{equation}{section}

\BibLatexMode{%
   \bibliography{template}
}

\title{Normality of closures of conormal bundles to Schubert cells for $\GL_n$}
\author{Ekaterina Bogdanova and Dyne Lee}
\date{}

\begin{document}

\maketitle

\begin{abstract}
    We prove normality of closures of the conormal bundle to certain Schubert cells in the flag variety for $\GL_n$ coming from Grassmannians. Our argument uses a Zelevinsky-type resolution of singularities, and through a sequence of cohomological computations based on the Borel–Weil–Bott theorem, we reduce normality to a problem about Kronecker coefficients, which we solve by applying results of Knutson and Tao from honeycomb theory.
\end{abstract}

\tableofcontents

\section{Introduction}
The goal of this paper is to study geometry of closures of conormal bundles to Schubert cells. These closures are the natural microlocal enhancements of the corresponding Schubert varieties, which has been extensively studied in {\cite{bri}, \cite{Sesh}, \cite{MeSri}, \cite{Raman}. The geometry of closures of conormal bundles to Schubert cells package information about singularities of Schubert varieties. Moreover, they are precisely the irreducible Lagrangian supports that occur in characteristic cycles of Schubert-constructible sheaves and $D$-modules, linking them to various topics in geometric representation theory. 

We prove that closures of conormal bundles to certain Schubert cells are normal in the case of general linear group, which is a particularly effective setting because the geometry admits explicit combinatorial and matrix descriptions. 

Let us formulate the main result precisely. Let $G=\GL_n$ be the general linear group over an algebraically closed field of characteristic 0 and $B\subseteq G$ the Borel subgroup of upper triangular matrices. Let $\mathcal{B}$ denote the flag variety. Recall that Schubert cells are indexed by elements of the Weyl group W, which is the symmetric group $S_n$ in our case. Denote the Schubert cell corresponding to $w \in W$ by $\mathcal{B}_w$.

\begin{notn}
    Let $s_i$ be the simple reflection transposing $i$ and $i+1$ and $w_0$ the longest element of the Weyl group ($w_0$ is the permutation which reverses the order of the elements $1,\ldots,n$). 
    
    For natural numbers $k$, $l$, $m$, with $k+l \leq n$, $m \leq \min(k, l)$, denote $\sigma_j=s_{k-j+1}s_{k-j+2}\ldots s_{n-l+m-j}$. We will call the elements of $W$ of the form
    \begin{align*}
    v_m=\sigma_1\sigma_2\ldots\sigma_mw_0
    \end{align*}
    \emph{good} elements of the Weyl group.
\end{notn}

Our main result is the following (Theorem \ref{main} in the main body of the paper):
\begin{thm*}
    For a good $w \in W$, the closure $\overline{T^*_{\mathcal{B}_w} \mathcal{B}}$ of the conormal to the Schubert cell corresponding to $w$ is normal.
\end{thm*}

We note that our techniques are different from those used in {\cite{bri}, \cite{Sesh}, \cite{MeSri}, \cite{Raman} to prove the normality of the closures of the Schubert cells themselves. 

\subsection{Idea of the proof} Fix a good element of $W$ and the corresponding numbers $k$, $l$ and $m$. First, we reduce the problem to the question of normality of the scheme $$\overline{T^*_{D_m}\Gr_k\times\Gr_l},$$ where $\Gr_i$ are Grassmannians parameterizing $i$-dimensional subspaces in the $n$-dimensional subspace, and $D_m$ is the $\GL_n$-orbit in $\Gr_k\times\Gr_l$ parameterizing pairs of subspaces with $m$-dimensional intersection. 

We prove normality of $\overline{T^*_{D_m}\Gr_k\times\Gr_l}$ using Stein factorization normality criterion. In other words, we construct a resolution of singularities $$p: \tilde{Z} \rightarrow \overline{T^*_{D_m}\Gr_k\times\Gr_l}$$
and prove that $p_* \mathcal{O}_{\tilde{Z}} \cong \mathcal{O}_{\overline{T^*_{D_m}\Gr_k\times\Gr_l}}$ (in fact, surjectivity suffices). Our resolution $\tilde{Z}$ imitates the resolution introduced by Zelevinsky in \cite{Zel}. In more detail, we show in Lemma \ref{closure} that 
\begin{equation}
    \overline{T^*_{D_m}\Gr_k\times\Gr_l} \cong \{(V_k,V_l,N)\in\Gr_k\times\Gr_l\times \overline O_m| \Imag(N)\subset V_k,V_l\subset\Ker(N),\dim(V_k\cap V_l)\geq m\},
\end{equation}
where $\overline O_m$ is the set of square-zero matrices of rank at most $m$.
Then $$\tilde{Z} := \{(L,V_k,V_l,H,N)\in\Gr_m\times\Gr_k\times\Gr_l\times\Gr_{n-m}\times\overline O_m: \Imag(N)\subset L\subset V_k,V_l\subset H\subset\Ker(N)\}.$$

Using Koszul resolutions, general derived pushforward criterions for comparing zeroth cohomology of non-derived pushforward and derived one, as well as Borel-Weil Bott theorem, we reduce the problem to purely combinatorial statement about Kronecker coefficients. Finally, we solve this problem using the results of Klyachko (\cite{Klya})and Knutson--Tao (\cite{KT}) relating non-vanishing of Littlewood-Richardson coefficient $c_{\nu, \tau}^{\mu}$ and existence of Hermitian matrices $A, B$ whose eigenvalues (in weakly decreasing order) are $nu$, $\tau$, and such that the eigenvalues of $A+B$ are $\mu$.

\subsection{Relationship to the work of Mehta--Trivedi.}\label{MT}
One can also get the result of the main theorem (Theorem \ref{main}) using \cite{MT}. We explain the reduction to Mehta--Trivedi's result in this subsection. However, we note that our proof is independent of their work and uses genuinely different methods.

Namely, in \cite[Theorem 4.1(3)]{MT} the authors prove that the ring $$\frac{K[S, N]}{(SN, NS, I_{m+1}(N), I_{\dim(A) -m+1}(S)}$$
is normal and Cohen-Macaulay. 

Let us show that Zariski-locally the scheme $$\overline{T^*_{D_m}\Gr_k\times\Gr_l}$$
is a product of an affine space and $$\Spec(\frac{K[S, N]}{(SN, NS, I_{m+1}(N), I_{\dim(A) -m+1}(S)}).$$ Indeed, note that $\overline{T^*_{D_m}\Gr_k\times\Gr_l}$ is equivalent to the subscheme of $\Gr_k\times\Gr_l\times \End(V)$ cut out by the following conditions on $(V_k, V_l, N)$:
    \begin{enumerate}
        \item $\Imag(N)\subset V_k,V_l\subset\Ker(N)$;
        \item all $(m+1)$-minors of $N$ vanish;
        \item all $(k+l-m+1)$-minors of the map 
        $$V_k \oplus V_l \rightarrow V$$
        sending $(v, k)$ to $v+k$ vanish.
    \end{enumerate}
    Choose a point $$(V_k^0, V_l^0, N^0) \in \overline{T^*_{D_m}\Gr_k\times\Gr_l}.$$ Consider the following Zariski open neighborhood of $(V_k^0, V_l^0, N^0)$. Consider a decomposition $$V \cong A \oplus B \oplus C \oplus D,$$ where $V_k^0 = A \oplus B$ and $V_l^0 = A \oplus C$. Intersect $$\overline{T^*_{D_m}\Gr_k\times\Gr_l}$$ with $$U \subset \Gr_k\times\Gr_l \times \End(V),$$ where 
     \[
    U=\left\{(V_k,V_l, N):
    \begin{array}{l}
    V=V_k\oplus(C\oplus D),\\
    V=V_l\oplus\bigl((V_k\cap(B\oplus C\oplus D))\oplus D\bigr)
    \end{array}
    \right\}.
    \]
   Note that the projection $H: V_k \rightarrow V/(C \oplus D)$ is an isomorphism. Hence $V_k = \Imag(\gra_H)$ for a unique $H \in \Hom(A \oplus B, C \oplus D)$. By definition, $\gra_H \oplus \Id_{C \oplus D}$ is invertible and 
   $$V = (\gra_H \oplus \Id_{C \oplus D})^{-1}(V_l) \oplus B \oplus D.$$
    Therefore $$(\gra_H \oplus \Id_{C \oplus D})^{-1}(V_l) = \Imag(\gra_S)$$
    for a unique $S \in \Hom(A \oplus C, B \oplus D)$. Without loss of generality (applying an automorphism preserving $A \oplus B$) we can assume that the image of $S$ projected to $B$ is $0$ and $S$ restricted to $C$ is $0$, i.e. $S: A \rightarrow D$.  Now the conditions defining $U$ imply that $N$ factors through $N: D \rightarrow A$. Now the conditions (1)-(3) defining $\overline{T^*_{D_m}\Gr_k\times\Gr_l}$ imply that 
    \begin{enumerate}
        \item $SN=0$;
        \item $NS=0$;
        \item $\rk(N) \leq m$;
        \item $\rk(S) \leq \dim(A) - m$.
    \end{enumerate}
    Therefore $$U \cap \overline{T^*_{D_m}\Gr_k\times\Gr_l}$$ is a product of an affine space and $$\Spec(\frac{K[S, N]}{(SN, NS, I_{m+1}(N), I_{\dim(A) -m+1}(S)}).$$

\subsection{Organization of the paper.} In section \ref{sec2}, we recall results on the cohomology of coherent sheaves on flag varieties and Grassmannians, which will be used in the proof of the main theorem. Also, we establish the moduli description of the closure to the conormal variety that we want to study. In section \ref{sec3}, we prove Theorem \ref{main}. We construct a resolution of singularities for the closure to the conormal variety, and reduce the problem to the question on Kronecker coefficients. In section \ref{sec4}, we establish the desired results about Kronecker coefficients.

\subsection{AI declaration.}
LLM's were used to get the results in Section \ref{sec4} of this text. LLM noticed the relationship of the present work to the work of \cite{MT}, explained in Subsection \ref{MT}. Everything else in the paper was produced by the authors.

LLM's were used to proof-read the earlier draft. This paper was written by the authors, and we claim responsibility for the correctedness of the presented results.

\subsection{Acknowledgements.} We want to thank Roman Bezrukavnikov for suggesting this problem and for many helpful conversations along the way. This project began as part of MIT SPUR program in the summer of 2025, and we would like to thank all organizers of this program. We thank Aleksandr Zimin for running his configuration of LLM's on Propositions \ref{kronecker} and \ref{kroneckergenweights}. We are grateful to Dora Woodruff for helpful discussions on Littlewood-Richardson and Kronecker coefficients.

E.B. thanks Lyalya Guseva for explaining the basics of coherent sheaves on Grassmannians.

\section{Preliminaries}\label{sec2}

\subsection{Notation}

Throughout this paper, we fix an algebraically closed field $K$ of characteristic 0 and a vector space $V$ of dimension $n$ over $K$. We denote by $G=\GL(V)=\GL_n$ the general linear group and $B\subset G$ a Borel subgroup. We identify the quotient $G/B$ with the flag variety $\mcB$ parametrizing complete flags of subspaces in $V$. 

$W=S_n$ is the Weyl group of $G$. Given an element $w\in W$, we denote by $C_w\subset\mcB$ the Schubert cell corresponding to $w$, which is a $B$-orbit in $\mcB$. We refer to \cite{bri} for further discussion on Schubert cells.

We denote by $\Gr(d, s)$ the Grassmannian of $d$-dimensional vector subspaces of $K^s$. When $s=n$, we write $\Gr_d$ for $\Gr(d, n)$. It is isomorphic to the homogeneous space $G/P_d$, where $P_d$ is the parabolic subgroup of block upper triangular matrices with two blocks of size $d$ and $n-d$. Since the Weyl group $W_{P_d}$ of $P_d$ is $S_d\times S_{n-d}$, $B$-orbits in $G/P_d$ correspond to elements of $S_n/(S_d\times S_{n-d})$ (See \cite{bri}).

\subsection{Schur functors}
Let \[
\lambda=(\lambda_1,\lambda_2,\ldots,\lambda_m)
\] be a nonincreasing sequence of nonnegative integers. We refer to $\lambda$ as a partition of $|\lambda|=\sum_i\lambda_i$, and we will often identify it with the Young diagram with rows of size $\lambda_i$. We denote by $\lambda^*$ the transpose of the Younge diagram.
Given an integer $a\in\Z$, we write $\lambda+a$ for the sequence \[
(\lambda_1+a,\lambda_2+a,\ldots,\lambda_m+a).
\] We denote by \[
\rho=(n-1,n-2,\ldots,0)
\] the sum of the fundamental weights of $G$.

The Schur functor $\Sigma^\lambda$ sends the vector space $V$ to the irreducible representation $\Sigma^\lambda V$ of $G=\GL(V)$ with highest weight $\lambda$. If $\lambda$ has more than $\dim V$ positive parts, $\Sigma^\lambda V=0$.
In particular, the $a$th symmetric power $S^aV$ is the Schur functor \[
\Sigma^{(a,0,\ldots,0)}V,
\] and the $a$th exterior power $\wedge^iV$ is \[
\Sigma^{(1,\ldots,1,0,\ldots,0)}V.
\] The character of the representation $\Sigma^\lambda V$ is the Schur polynomial $s_\lambda$.

We write $-\lambda$ for the nonincreasing sequence \[
(-\lambda_m,\ldots,-\lambda_2,-\lambda_1)
\] and \[
\Sigma^{-\lambda}V=\Sigma^\lambda V^\vee.
\]
We refer to \cite[Chapter~6]{fh} for further discussion on Schur functors.

Let $\mcV$ be a vector bundle of rank $r$ on a scheme $X$. We may apply the Schur functor $\Sigma^\lambda$ to $\mcV$ to obtain the vector bundle corresponding to the representation of $\GL_r$ with highest weight $\lambda$.

\subsection{Coherent sheaves on flag varieties}

On the Grassmannian $\Gr_d$, we denote by $\underline V$ the trivial vector bundle corresponding to the vector space $V$, $\mcS$ the tautological vector bundle of rank $d$, and $\mathcal Q$ the quotient bundle $\underline V/\mcS$.

Given an integer sequence $\alpha$ of length $n$, there is a one-dimensional representation of $G$ given by the character $\alpha$ of $B$. Thus on the complete flag variety $G/B$, we obtain a line bundle $\mcO(\alpha)$. The cohomology of $\mcO(\alpha)$ is computed as follows (see \cite{dem}):

There is at most one permutation $\sigma\in S_n$ such that $\sigma(\alpha+\rho)-\rho$ is a nonincreasing integer sequence (i.e. it is a dominant weight of $G$). Such a $\sigma$ exists if and only if $\alpha+\rho$ has no repeated parts, and in this case $\sigma$ is the permutation putting $\alpha+\rho$ into strictly decreasing order.

\begin{thm}[Borel-Weil-Bott Theorem] \label{bwb}
    (a) If there is no $\sigma\in S_n$ such that $\sigma(\alpha+\rho)-\rho$ is nonincreasing, $\mcO(\alpha)$ has no cohomology in any degree.

    (b) If such a $\sigma$ exists, $\mcO(\alpha)$ has cohomology only in degree equal to the length $l(\sigma)$, and \[
    H^{l(\sigma)}(G/B,\mcO(\alpha))=\Sigma^{\sigma(\alpha+\rho)-\rho}V^\vee.
    \]
\end{thm}

    Let $Q$ be the incidence subvariety in the product $\Gr_k \times \Gr_l$. 
\begin{cor}
    Consider the vector bundle $\Sigma^{\alpha} \mathcal{S}_k^\vee \otimes \Sigma^{\beta} \mathcal{Q}_l^\vee.$ Then we compute its cohomology in the following way:
    \begin{enumerate}
        \item If there is no $\sigma\in S_n$ such that $\sigma((\alpha, \beta)+\rho)-\rho$ is nonincreasing, $\mcO(\alpha)$ has no cohomology in any degree.
        \item If such a $\sigma$ exists, $\Sigma^{\alpha} \mathcal{S}_k^\vee \otimes \Sigma^{\beta} \mathcal{Q}_l^\vee$ has cohomology only in degree equal to the length $l(\sigma)$, and \[
    H^{l(\sigma)}(Q,\Sigma^{\alpha} \mathcal{S}_k^\vee \otimes \Sigma^{\beta} \mathcal{Q}_l^\vee)=\Sigma^{\sigma(\alpha+\rho)-\rho}V^\vee.
    \]
    \end{enumerate}
\end{cor}

\begin{lem}\label{bwbcor}
    For Young diagrams $\lambda$ and $\mu$, The vector bundle $$\Sigma^{\lambda} \mathcal{S}_m \otimes \Sigma^{\mu} \mathcal{Q}_{n-m}^\vee \otimes \Sym^j(\mathcal{S}_m^\vee \otimes \mathcal{Q}_{n-m})$$ on the incidence subvariety of $\Gr_m \times \Gr_{n-m}$ has cohomology in degree at most $|\lambda|+|\mu|-1$.
\end{lem}
\begin{proof}
    Recall that  \[
     \Sym^j(V\otimes W)=\bigoplus_{\lambda\vdash j}\Sigma^\lambda V\otimes\Sigma^\lambda W.
     \] 
     By Pieri's formula,  each irreducible summand of the term \[
     \Sigma^{\lambda}\mcS_m\otimes\Sigma^\mu \mcQ_{n-m}^\vee\otimes\Sigma^{\alpha}\mcS_m^\vee\otimes\Sigma^\alpha \mcQ_{n-m}
     \] 
     has the form
     \begin{equation}\label{schur}
         \Sigma^{\alpha - (a_1, \ldots , a_m)}\mcS_m^\vee \otimes \Sigma^{\alpha - (b_1, \ldots b_m} \mcQ_{n-m}
     \end{equation}
     with $a_i, b_i \geq 0$ and $\Sigma a_i = | \lambda |$ and $\Sigma b_i = | \mu |$.

     Theorem \ref{bwb} implies that cohomology of (\ref{schur}) may occur only in degree equal to the number of inversions in the sequence 
     \begin{equation}
        (\alpha_1, \ldots, \alpha_m, 0, \ldots, 0, -\alpha_m, \ldots, -\alpha_1)  + (-a_1, \ldots, -a_m, 0, \ldots, 0, b_1, \ldots, b_m). 
     \end{equation}
     However, we claim that for a sequence
$$x_1>\cdots>x_n$$
and $y_1,\ldots,y_n$ distinct integers such that $y\ne x$ we have the following inequality of the number of inversions:
$$
\operatorname{inv}(y)<\sum_t|y_t-x_t|.$$
The result follows.
\end{proof}

\subsection{The conormal bundle}
\begin{lem}\label{conormal}
    (1) The cotangent space $T^*(\Gr_k\times\Gr_l)$ at the point $(V_k,V_l)$ is identified with \begin{align*}
        \Hom(V/V_k,V_k)\oplus\Hom(V/V_l,V_l).
    \end{align*}

    (2) The fiber of the conormal sheaf $T^*_{D_m}(\Gr_k\times\Gr_l)$ at the point $(V_k,V_l)\in D_m$ is identified with \begin{align*}
        \Hom(V/(V_k+V_l),V_k\cap V_l).
    \end{align*}
\end{lem}
 
\begin{proof}
    (1) The tangent space to $\Gr_k$ at a point identified with the subspace $V_k$ is $\Hom(V_k,V/V_k)$ (see \cite{gw}). Thus the cotangent space is $\Hom(V/V_k,V_k)$, and the statement follows. \\

    (2) Fix a point $x=(V_k,V_l)\in D_m$. There is a short exact sequence \begin{align*}
        0\rightarrow T(D_m)_x\rightarrow T(\Gr_k\times\Gr_l)_x\rightarrow T_{D_m}(\Gr_k\times\Gr_l)_x\rightarrow0.
    \end{align*} We compute $T(D_m)_x$. Consider the tangent space as a set of $K[\epsilon]/\epsilon^2$-valued points of $D_m$. A point in this tangent space is identified with a pair of projective submodules $\mcV_k,\mcV_l\subset V\otimes K[\epsilon]$ of dimensions $k$ and $l$ satisfying the following conditions: 
    \begin{enumerate}
        \item $\mcV_k \cap \mcV_l$ is projective of dimension $m$, 
        \item $V\otimes K[\epsilon]/\mcV_k$, $V\otimes K[\epsilon]/\mcV_l$, and $V\otimes K[\epsilon]/(\mcV_k \cap \mcV_l)$ are projective,
        \item $\mcV_k/\epsilon\simeq V_k$  and $\mcV_l/\epsilon\simeq V_l$, and (4) the map \begin{align*}
        \mcV_k\cap\mcV_l\rightarrow\mcV_k/\epsilon\cap \mcV_l/\epsilon
    \end{align*} is surjective. 
    \end{enumerate}
    
    We note that $\mcV_k$ is the image of a $K[\epsilon]$-linear map $V_k\otimes K[\epsilon]\rightarrow V\otimes K[\epsilon]=V\oplus \epsilon V$, and likewise for $\mcV_l$. Given a map $w:V_k\rightarrow V$, there is a unique map \begin{align*}
        \tilde w: V_k\otimes K[\epsilon]\rightarrow V\oplus \epsilon V
    \end{align*} such that $\tilde w|_{V_k}=\iota\oplus\epsilon w$, where $\iota:V_k\rightarrow V$ is the inclusion. Then $\Imag(\tilde w)$ is some projective module $\mcV_k$. The same statement holds for $V_l$.
    Thus there is a map \begin{align*}
        \Hom(V_k,V)\oplus\Hom(V_l,V)&\rightarrow T(\Gr_k\times\Gr_l)_x,\\
        (w_k,w_l)&\mapsto(\Imag(\tilde w_k),\Imag(\tilde w_l)).
    \end{align*}
    Its kernel consists of the pairs of maps $(w_k,w_l)$ such that $w_k(V_k)\subset V_k$ and $w_l(V_l)\subset V_l$. Imposing the surjectivity condition on $\mcV_k\cap\mcV_l$, the tangent space $T(D_m)_x$ is identified with pairs $(f,f')\in\Hom(V_k,V/V_k)\oplus\Hom(V_l,V/V_l)$ which agree on $V_k\cap V_l$ after projection to $V/V_k\cap V/V_l=V/(V_k+V_l).$ Now the quotient of $T(\Gr_k\times\Gr_l)_x$ by $T(D_m)_x$ can be identified with \begin{align*}
        \Hom(V_k\cap V_l,V/(V_k+V_l)).
    \end{align*} The dual of this space is the fiber of the conormal bundle at $x$.
\end{proof}

We now give a useful description of the closure of the conormal to an orbit inside the tangent bundle. Define  \[
O_m=\{N\in\Mat_n|N^2=0,\rk N=m\}.
\] Let $r = \min (m, \lfloor \frac{n}{2} \rfloor).$ Note that $O_m$ is the nilpotent orbit of Jordan type $(2^r, 1^{n-2r})$. Its closure $\overline O_m$ is the set \[
\overline O_m=\{N\in\Mat_n|N^2=0,\rk N\leq m\}.
\]

\begin{lem}\label{closure}
    Let $Y_m=T^*_{D_m}(\Gr_k\times\Gr_l)$ be the conormal bundle to the $G$-orbit $D_m$ in $\Gr_k\times\Gr_l$. Its closure $\overline Y_m\subset T^*(\Gr_k\times\Gr_l)$ is isomorphic to the reduced subvariety of \[
    Z=\{(V_k,V_l,N)\in\Gr_k\times\Gr_l\times \overline O_m| \Imag(N)\subset V_k,V_l\subset\Ker(N),\dim(V_k\cap V_l)\geq m\}.
    \]
\end{lem}

\begin{rem}
    By definition, $k$ and $l$ are always greater than or equal to $m$, and we may observe that they are also at most $n-m$. Thus $Z$ is nonempty.
\end{rem}

\begin{proof}
    Given a point $(V_k,V_l,N)\in Z$, we view $N$ as a map $V/V_k\to V_k$ and also as a map $V/V_l\to V_l$.
    Then viewing the cotangent space as in Lemma \ref{conormal}, we embed $Y_m$ into $T^*(\Gr_k\times\Gr_l)$ so that $(V_k,V_l,N)$ is sent to the $(N,-N)$ in the cotangent space at $(V_k,V_l)$. Now $Y_m$ consists of the triples $(V_k,V_l,N)$ where $\Imag(N)\subset V_k,V_l\subset\Ker(N)$, $V_k\cap V_l$ has dimension $m$, and $N$ is nilpotent of rank at most $m$, again by Lemma \ref{conormal}. 
    
    We show that the closure of the above space is $Z$. Consider the projections from $Y_m$ and $Z$ to $\overline{O}_m$, each forgetting everything but $N$. Given $N\in {O}_m$, the fiber over $N$ in $Y_m$ is the set of pairs of subspaces $(W_{k-m},W_{l-m})$ in $\Ker(N)/\Imag(N)$ of rank $k-m$ and $l-m$, respectively, such that $W_{k-m}\cap W_{l-m}=0$. The fiber over $N$ in $Z$ is the same set of pairs without the condition that the subspaces do not intersect. This is an open non-empty condition, hence it is dense since the product of Grassmannians is irreducible.

    It remains to show that $Z$ is irreducible. We claim that the open locus $Z^\circ \subset Z$ where $N \in O_m$ is dense in $Z$. Indeed, for a fixed pair $(V_k,V_l)$ the fiber $Z_{(V_k,V_l)}^\circ \subset Z_{(V_k,V_l)}$ is open and dense.
\end{proof}

\section{Construction of the Resolution}\label{sec3}

\subsection{Reduction to $\Gr_k\times\Gr_l$}\label{reduction}
In this subsection we reduce the question of normality of the closure $\overline{T^*_{\mathcal{B}_w} \mathcal{B}}$ for a good $w \in W$ to the question of normality of conormal closure to some $G$-orbit in the product of two Grassmannians.

Consider the diagonal action of $G$ on $\mcB^2$. We describe $G$-orbits under this action. Given any point $(F_1,F_2)$ in $\mcB^2$, there exists $g\in G$ such that $g(F_1,F_2)=(F_s,F')$ where $F_s$ is the standard flag. Then $B$ preserves the first factor, so in the second factor we may obtain any flag in the $B$-orbit of $F\in\mcB$. Hence the $G$-orbits in $\mcB^2$ are also labeled by elements of $W$. 

Now consider the map $\mcB^2\to\mcB$ which is projection on the first factor of $\mcB$. The fiber over the standard flag $F_s$ may be identified with $\mcB$, so there is an inclusion $\iota: \mcB\to\mcB^2$. Given an element $w\in W$, the Schubert cell $C_w$ is the preimage of the $G$-orbit in $\mcB^2$ corresponding to $w$. Thus to show that the closure of $T^*_{C_w}\mcB$ is normal, it suffices to show that the closure of conormal to the corresponding $G$-orbit is normal.
    
    Suppose $k+l\leq n$, and let $\pr:\mcB^2\rightarrow\Gr_k\times\Gr_l$ be the map forgetting everything except the $k$-dimensional subspace of the first flag and the $l$-dimensional subspace of the second flag. Through the diagonal action of $G$, we may always standardize the first subspace to be $\left<e_1,\ldots,e_k\right>$. Thus $G$-orbits in $\Gr_k\times\Gr_l$ are parametrized by $W_{P_k}\backslash W/W_{P_l}$. These double cosets are of the form $X_i$ for $0\leq i\leq\min(k,l)$, where $X_i$ contains the permutations where $i$ of the first $l$ elements are sent to the first $k$ positions. This corresponds to the $G$-orbit of pairs of subspaces with exactly $i$-dimensional intersection.

     The preimage of $\overline{D_m}$ along $\pr$ is the union of some $G$-orbits in $\mcB^2$, one of which is open and dense in $\pr^{-1}(\overline{D_m})$. This open orbit corresponds to the element of the Weyl group which is last in the Bruhat order among these orbits. We denote this element by $v_m$. Then since the closure of $\mcB^{2,0}_{v_m}$ is $\pr^{-1}(D_m)$, the closures of the conormals to $\mcB^{2,0}_{v_m}$ and $\pr^{-1}(D_m)$ are isomorphic. Thus it suffices to show that $\overline{T^*_{\pr^{-1}(D_m)}\mcB^2}$ is normal, which occurs when \begin{align*}
    \overline{T^*_{D_m}\Gr_k\times\Gr_l}
\end{align*} is normal. 

An explicit description of $v_m$ is as follows:
Let $s_i$ be the simple reflection transposing $i$ and $i+1$, and $w_0$ the longest element of the Weyl group, which reverses the order of the elements $1,\ldots,n$. We denote $\sigma_j=s_{k-j+1}s_{k-j+2}\ldots s_{n-l+m-j}$. Then \begin{align*}
    v_m=\sigma_1\sigma_2\ldots\sigma_mw_0.
\end{align*}

The rest of the section is dedicated to the proof of normality of $Y_m$.

\subsection{Cohomological vanishing} In this subsection we perform some spectral sequence calculations that will be used in the proof of the main theorem.

\begin{lem}\label{spectral1}
    Let $f:X\to Y$ be a morphism of ringed spaces. Given a complex of $\mcO_X$-modules \[
    0\to M^{-n}\to\ldots\to M^{-1}\to M^0\to 0,
    \] if \[
    R^{i}f_*(M^{-i})=0
    \] for all $i>0$, then \[
    H^0(f_*M^{\bullet})\to H^0(Rf_*M^\bullet)
    \] is surjective.
\end{lem}

\begin{proof}
    Let $E$ be the hypercohomology spectral sequence such that \[
    E^{p,q}_1=R^qf_*M^p,
    \]
    which converges to $H^{p+q}(Rf_*M^\bullet)$. The condition that $R^if_*(M^{-i})=0$ implies that $E^{-i,i}_1=0$ and thus $E^{-i,i}_\infty=0$ for all $i>0$. 
    Therefore $E_2^{0,0} \cong H^0(f_* M^\bullet)$ surjects onto $E_{\infty}^{0,0} = H^0(Rf_* M^{\bullet})$.
\end{proof}

\begin{lem}\label{spectral2}
    Let $f:X\to Y$ be a morphism of ringed spaces. Given a complex of $\mcO_X$-modules \[
    0\to M^{-n}\to\ldots\to M^{-1}\to M^0\to 0,
    \] if \[
    R^{i+1}f_*(H^{-i}(M^\bullet))=0
    \] for all $i>0$, then \[
    H^0(Rf_*M^\bullet)\rightarrow f_*(H^0(M^\bullet))
    \] is surjective.
\end{lem}

\begin{proof}
    Let $E$ be the hypercohomology spectral sequence with \[
    E_2^{p,q}=R^qf_*H^p(M^\bullet),
    \] which converges to $H^{p+q}(Rf_*M^\bullet)$. There are maps \begin{align*}
        E^{-2,1}_2\rightarrow E^{0,0}_2\rightarrow E^{2,-1}_2.
    \end{align*} The leftmost term is always 0, and the condition that $R^{i+1}f_*(H^{-i}(M^\bullet))=0$ implies that the rightmost term is 0. Repeating the above on every page, the condition ensures that $E^{0,0}_2=E^{0,0}_\infty$. Now there is a filtration \[
    0=F^sH^0(Rf_*M^{\bullet})\subseteq F^{s-1}H^0(Rf_*M^{\bullet})\subseteq\ldots\subseteq 
H^0(Rf_*M^{\bullet})
    \] such that $E^{0,0}_\infty=F^0H^0(Rf_*M^{\bullet})/F^1H^0(Rf_*M^{\bullet})$. Since $E^{p,-p}_2$ and thus $E^{p,-p}_\infty$ are 0 for all $p<0$, $F^pH^0(Rf_*M^{\bullet})=0$ for all $p<0$. Thus $F^0H^0(Rf_*M^{\bullet})=H^0(Rf_*M^{\bullet})$, so $H^0(Rf_*M^{\bullet})$ surjects onto $E^{0,0}_\infty=f_*(H^0(M^\bullet))$.
\end{proof}

\subsection{Proof of the main result}
For natural numbers $k$, $l$, $m$, with $k+l \leq n$, $m \leq \min(k, l)$, denote $\sigma_j=s_{k-j+1}s_{k-j+2}\ldots s_{n-l+m-j}$. Let
    \begin{align*}
    v_m=\sigma_1\sigma_2\ldots\sigma_mw_0
    \end{align*}

\begin{thm}\label{main}
The closure of the conormal bundle to the Schubert cell $\mathcal{B}_{v_m}$ $$\overline{T^*_{\mathcal{B}_w} \mathcal{B}}$$ is a normal variety.
\end{thm}

The discussion in Section \ref{reduction} reduces this statement to showing that
$\overline{T^*_{D_m}(\Gr_k\times\Gr_l)}$ is normal. By
Lemma \ref{closure}, this variety is identified, with its reduced structure, with
\[
Z=\{(V_k,V_l,N)\in \Gr_k\times\Gr_l\times\overline O_m :
\Imag N\subset V_k\cap V_l,\ V_k+V_l\subset\Ker N,\
\dim(V_k\cap V_l)\ge m\}.
\]

We begin by constructing a resolution of singularities for $Z$. Throughout, we use $m\le\min(k,l)$ and $k+l\le n$. In particular,
$2m\le k+l\le n$, so $m\le n-m$, and $k,l\le n-m$.

\begin{prop}\label{resolution}
    Consider
\[
\tilde Z=\{(L,V_k,V_l,H,N)\in \Gr_m\times\Gr_k\times\Gr_l
\times\Gr_{n-m}\times\overline O_m :
\Imag N\subset L\subset V_k\cap V_l,\ V_k+V_l\subset H\subset\Ker N\},
\]
and let $p:\tilde Z\to Z$ be the map forgetting $L$ and $H$.
Then  $p$ is a resolution of singularities, i.e. $\tilde Z$
is smooth and irreducible and $p$ is projective and birational.
\end{prop}

\begin{proof}
    We note that the morphism $p$ is well-defined since 
if $(L,V_k,V_l,H,N)\in\tilde Z$, then $\dim(V_k\cap V_l)\geq\dim L=m$.
    Define
\[
X=\{(L,V_k,V_l,H)\in\Gr_m\times\Gr_k\times\Gr_l\times\Gr_{n-m}
: L\subset V_k\cap V_l,\ V_k+V_l\subset H\}.
\] We first show that $X$ is smooth, irreducible and projective.
The projection $X\to Q=\{(L,H): L\subset H\}$ is $G$-equivariant. Here $Q$ is
a partial flag variety (nonempty since $m\le n-m$), hence smooth, projective
and irreducible. The fiber over $(L,H)$ is
\[
\Gr(k-m,H/L)\times\Gr(l-m,H/L),\qquad \dim H/L=n-2m,
\]
which is nonempty because $k-m,\,l-m\le n-2m$. So $X\to Q$ is a Zariski
locally trivial fiber bundle (it is $G\times^{P}F$ for the stabilizer $P$ of a
point of $Q$) with smooth, irreducible, projective fiber. Hence $X$ is smooth,
irreducible and projective.

Now we show that $\tilde Z$ is smooth and irreducible.
Given $(L,V_k,V_l,H)\in X$, the endomorphisms $N$ of $V$ with
$\Imag N\subset L$ and $H\subset\Ker N$ are exactly the
compositions $V\to V/H\to L\hookrightarrow V$. Every such $N$ satisfies
$N^2=0$ (since $L\subset H$) and $\rk N\le m$, so it
automatically lies in $\overline O_m$. Therefore $\tilde Z$ is
the total space of the vector bundle $\SHom(V/H,L)$ over $X$, and in
particular it is smooth and irreducible.

Furthermore, $p$ is projective and surjective.
The conditions defining $\tilde Z$ are closed, so $\tilde Z$ is a closed
subvariety of $\Gr_m\times\Gr_{n-m}\times Z$. Thus $p$ is the
restriction of a projection with projective fibers, and is projective. For
surjectivity, let $(V_k,V_l,N)\in Z$.
\begin{itemize}
\item Since $\rk N\le m\le\dim(V_k\cap V_l)$ and
$\Imag N\subset V_k\cap V_l$, there is an $m$-plane $L$ with
$\Imag N\subset L\subset V_k\cap V_l$.
\item Since
\[
\dim(V_k+V_l)=k+l-\dim(V_k\cap V_l)\le k+l-m\le n-m\le n-\rk N=\dim\Ker N,
\]
there is an $(n-m)$-plane $H$ with $V_k+V_l\subset H\subset\Ker N$.
\end{itemize}
Then $(L,V_k,V_l,H,N)\in\tilde Z$ maps to $(V_k,V_l,N)$.

We now show that $p$ is an isomorphism over the rank-$m$ locus.
Let $Z^\circ=\{(V_k,V_l,N)\in Z:\rk N=m\}$, which is open in $Z$.
Over $Z^\circ$ the conditions on $\tilde Z$ force
$L=\Imag N$ and $H=\Ker N$ for dimension reasons.
Because $N$ has constant rank on $Z^\circ$, the assignments
$N\mapsto\Imag N$ and $N\mapsto\Ker N$ define morphisms
$Z^\circ\to\Gr_m$ and $Z^\circ\to\Gr_{n-m}$. Hence
\[
(V_k,V_l,N)\mapsto(\Imag N,V_k,V_l,\Ker N,N)
\]
is a morphism inverse to $p$, and $p$ restricts to an isomorphism
$p^{-1}(Z^\circ)\xrightarrow{\ \sim\ }Z^\circ$.

We note that the locus $\{\dim(V_k\cap V_l)=m\}$ is not in general a
locus where $p$ is bijective. When $k+l<n$ and $\rk N<m$, the
plane $H$ can usually be chosen in more than one way.

In the proof of Lemma \ref{closure} we showed that $Z^\circ$ is dense in $Z$, so $p$ is birational.
\end{proof}

Now we proceed with the proof of the main theorem.
\begin{proof} (Theorem \ref{main})
We have shown that $p: \tilde Z\to Z$ is a projective birational
morphism from a smooth irreducible variety, i.e. a resolution of
singularities. Consider its Stein factorization
\[
\tilde Z\to\operatorname{\mathbf{Spec}}_Z p_*\mathcal O_{\tilde Z}\to Z_{\operatorname{red}}.
\]
The second map is finite and birational, and its source is normal because
$\tilde Z$ is normal. So it is the normalization of $Z_{\operatorname{red}}$. Consequently $Z_{\operatorname{red}}$ is
normal if and only if $$\mathcal O_{Z_{\operatorname{red}}}\to p_*\mathcal O_{\tilde Z}$$ is an
isomorphism. This map is injective because by Proposition $Z_{\operatorname{red}}$ is reduced and $p$ is
surjective, so it suffices to prove that it is surjective. It suffices to prove that $$\mathcal{O}_Z \rightarrow p_*\mathcal O_{\tilde Z}$$ is surjective.

We also define \begin{align*}
    \tilde Q&=\{(L,H,N)\in\Gr_m\times\Gr_{n-m}\times\overline O_m: \Imag(N)\subset L\subset H\subset\Ker(N)\},\\
        Q&=\{(L,H)\in\Gr_m\times\Gr_{n-m}: L\subset H\}.
\end{align*} Now there is a commutative diagram \begin{equation}\label{diagram}
    \begin{tikzcd}
        &X\arrow[hookrightarrow,r,"h"]&Q\times\Gr_k\times\Gr_l\arrow[rd]\\
        &\tilde Z\arrow[u]\arrow[hookrightarrow,r,"i"]\arrow[d,"p"]&\tilde Q\times\Gr_k\times\Gr_l\arrow[u,"\phi"]\arrow[d,"\varpi"]&\Gr_k\times\Gr_l\\
        &Z\arrow[hookrightarrow, r,"j"]&\overline O_m\times\Gr_k\times\Gr_l\arrow[ru]
    \end{tikzcd}
    \end{equation} where the upper square is Cartesian.

     Since $j$ in the diagram \ref{diagram} is a closed embedding, it suffices to show that \[
    j_*\mcO_Z \to j_*p_*\mcO_{\tilde Z}
    \] is surjective. 
    
    We note that $\overline O_m\times\Gr_k\times\Gr_l$ is normal (see \cite{kp} for the normality of $\overline O_m$) and that $\varpi$ is a proper birational morphism with connected nonempty fibers. Hence its Stein factorization implies that \[
    \varpi_*\mcO_{\tilde Q\times\Gr_k\times\Gr_l} = \mcO_{\overline O_m\times\Gr_k\times\Gr_l}.
    \] Then using the lower square in the diagram, it suffices to show that \begin{align*}
        \varpi_*\mcO_{\tilde Q\times\Gr_k\times\Gr_l}\to \varpi_*i_*\mcO_{\tilde Z}
    \end{align*} is surjective. We proceed by constructing a complex of sheaves on $\tilde Q\times\Gr_k\times\Gr_l$ ending with the map \begin{align*}
        \mcO_{\tilde Q\times\Gr_k\times\Gr_l}\rightarrow i_*\mcO_{\tilde Z}
    \end{align*} and considering the pushforward of this complex along $\varpi$.

    Define the subvariety $X_k\subset Q\times\Gr_k$ of triples $(L,H,V_k)$ such that $L\subset V_k\subset H$. We have a projection $p_{m,k}: Q\times\Gr_k\rightarrow\Gr_m\times\Gr_k$, where the image of $X_k$ is the set $S_{m,k}$ of pairs $(L,V_k)$ with $L\subset V_k$. There is a regular section of the vector bundle $E_{k,1}=\SHom(L,V/V_k)$ on $\Gr_m\times\Gr_k$ whose zero locus is $S$. The rank of this bundle is $r_1=m(n-k)$, which is equal to the codimension of $S$ in $\Gr_m\times\Gr_k$. Thus we have a Koszul resolution \begin{align}\label{koszul1}
        0\rightarrow\wedge^{r_1}E_{k,1}^\vee\rightarrow\ldots \rightarrow E_{k,1}^\vee\rightarrow\mcO_{\Gr_m\times\Gr_k} \rightarrow0,
    \end{align} where the rightmost cohomology $H^0$ is $\mcO_{S_{m,k}}$ and the complex is exact everywhere else.

    Similarly, we can project $p_{n-m,k}:Q\times\Gr_k\to\Gr_{n-m}\times\Gr_k$ so that the image of $X_k$ is $S_{n-m,k}$. Note that $X_k=Y_1\cap Y_2$, where \begin{align}\label{ydef}
        Y_1&=p_{m,k}^{-1}(S_{m,k})\\
        Y_2&=p_{n-m,k}^{-1}(S_{n-m,k}).
    \end{align}
    we obtain a Koszul resolution for the set of $(V_k,H)$ with $V_k\subset H$ in $\Gr_k\times\Gr_{n-m}$ using the vector bundle $E_{k,2}=\SHom(V_k,V/H)$. Tensoring the two resolutions for $E_{k,1}$ and $E_{k,2}$, we obtain a complex $M_k^\bullet$ with $H^0(M_k^\bullet)=\mcO_{X_k}$. However, since this is a Koszul complex of a vector bundle $E_{k,1}\oplus E_{k,2}$ of rank $nm$ and the codimension of $X_k$ in $Q\times\Gr_k$ is $m(n-m)$, the complex also has cohomology in degrees up to $-m^2$. 

    We do the same after replacing $\Gr_k$ with $\Gr_l$ to obtain a complex $M_l^\bullet$. Tensoring $M_k^\bullet$ and $M_l^\bullet$ yields a complex \begin{align*}
        M^\bullet=0\rightarrow\wedge^r E^\vee\rightarrow\ldots\rightarrow E^\vee \rightarrow\mcO_{Q\times\Gr_k\times\Gr_l}\to h_*\mcO_X
        \rightarrow0.
    \end{align*} Pulling back $M^\bullet$ along $\phi$ results in a complex of sheaves on $\tilde Q\times\Gr_k\times\Gr_l$, with the cohomology at degree 0 being $j_*\mcO_{\tilde Z}$. It now suffices to show that further taking the pushforward of this complex along $\varpi$ preserves $H^0$.

    Suppose that the conditions of Lemma \ref{spectral1} and Lemma \ref{spectral2} hold for $\phi^*M^\bullet$ and $\varpi$. Then we have that \begin{align*}
        H^0(\varpi_*\phi^*M^\bullet)\rightarrow\varpi_*(H^0(\phi^*M^\bullet))
    \end{align*} is surjective. Since the pushforward of the structure sheaf along $\varpi$ is the structure sheaf on $$\overline O_m\times\Gr_k\times\Gr_l,$$ this yields our desired surjectivity. Thus Propositions \ref{vanishing1} and \ref{vanishing2}, which show that the conditions of Lemma \ref{spectral1} and Lemma \ref{spectral2} are satisfied, finish the proof.
\end{proof}

\begin{prop}\label{vanishing1}
     Let $M^\bullet$, $\varpi$, and $\phi$ be defined as above. \begin{align*}
        R^{i}\varpi_*(\phi^*M^{-i})&=0
    \end{align*} for all $i>0$.
\end{prop}

\begin{proof}
    The complex $M^\bullet$ is obtained by taking the tensor product of four Koszul complexes as described in the proof of Theorem \ref{main}. $\phi^*M^{-i}$ is the direct sum of terms of the form \begin{align*}
        \mcF=\wedge^{i_1}(\mcS_m\boxtimes \mcQ_k^\vee)\otimes \wedge^{i_2}(\mcS_k\boxtimes \mcQ_{n-m}^\vee)\otimes \wedge^{i_3}(\mcS_m\boxtimes \mcQ_l^\vee)\otimes \wedge^{i_4}(\mcS_l\boxtimes \mcQ_{n-m}^\vee),
    \end{align*} where $\mcS_j$ and $\mcQ_j$ denote the tautological bundle and the tautological quotient bundle, respectively, on $\Gr_j$, and $i_1+i_2+i_3+i_4=i$. We rewrite this sheaf in terms of Schur functors, using the fact (see \cite[Exercise~6.11]{fh}) that \begin{align*}
        \wedge^i(A\otimes B)=\bigoplus_{|\lambda|=i}\Sigma^\lambda A\otimes\Sigma^{\lambda^*}B.
    \end{align*} 
    Thus it suffices to show that \begin{align*}
        R^{i}\varpi_*(\Sigma^{\lambda_1}\mcS_m\otimes\Sigma^{\lambda_3}\mcS_m\otimes \Sigma^{\lambda_2}\mcQ_{n-m}^\vee\otimes\Sigma^{\lambda_4}\mcQ_{n-m}^\vee\otimes F)=0
    \end{align*} for all $|\lambda_j|=i_j$, where $F$ denotes any factors of $\mcS_k,\mcS_l,\mcQ_k^\vee,\mcQ_l^\vee$ coming from the wedge power expansion.
    The decomposition of tensor products of Schur functors is described by the Littlewood-Richardson rule. We will work with the rule in terms of expansion of Young diagrams, described in \cite[Appendix~A.1]{fh}. We decompose a tensor product \[\Sigma^\lambda\mcF\otimes\Sigma^\mu\mcF=\bigoplus_{\alpha}\Sigma^{\alpha}\mcF
    \] for Young diagrams $\alpha$ obtained by adding the boxes of $|\mu|$ to $|\lambda|$ under certain combinatorial conditions. Here, we need only the condition that \[
    |\alpha|=|\lambda|+|\mu|.
    \]
    
    Now $\phi^*M^{-i}$ is a finite direct sum of bundles of the form
    \[
    \mathcal A_{\lambda,\mu}\boxtimes\mathcal B,
    \qquad
    \mathcal A_{\lambda,\mu}=f^*\big(\Sigma^\lambda\mathcal S_m\otimes\Sigma^\mu\mathcal Q_{n-m}^\vee\big),
    \qquad |\lambda|+|\mu|=i.
    \]
    Here $f:\tilde Q\to Q$ is the projection, and $\mathcal B$ is a vector
    bundle on $\Gr_k\times\Gr_l$ built from Schur functors of
    $\mathcal S_k,\mathcal S_l,\mathcal Q_k^\vee,\mathcal Q_l^\vee$. Since
    $R^i\varpi_*$ commutes with finite direct sums, it suffices to show that
    \[
    R^i\varpi_*\big(\mathcal A_{\lambda,\mu}\boxtimes\mathcal B\big)=0
    \qquad\text{whenever } |\lambda|+|\mu|=i>0.
    \]
    
    Write $\psi:\tilde Q\to\overline{O}_m$ for the map
    $(L,H,N)\mapsto N$, so that
    $\varpi=\psi\times\mathrm{id}_{\Gr_k\times\Gr_l}$. The
    projection $\overline{O}_m\times\Gr_k\times\Gr_l\to\overline{O}_m$
    is flat, so flat base change gives
    $R\varpi_*(\mathrm{pr}_1^*\mathcal A)=\mathrm{pr}_1^*R\psi_*\mathcal A$.
    Since $\mathcal B$ is locally free, the projection formula then gives
    \[
    R^i\varpi_*\big(\mathcal A_{\lambda,\mu}\boxtimes\mathcal B\big)
    \cong R^i\psi_*\mathcal A_{\lambda,\mu}\boxtimes\mathcal B .
    \]
    It is therefore enough to show that $R^i\psi_*\mathcal A_{\lambda,\mu}=0$.
    
    Since $\overline{O}_m$ is affine, $R^i\psi_*\mathcal A_{\lambda,\mu}=0$ if and only if
    \[
    H^i\big(\tilde Q,\,f^*(\Sigma^\lambda\mathcal S_m\otimes\Sigma^\mu\mathcal Q_{n-m}^\vee)\big)=0.
    \]
    
    Next, $\tilde Q$ is the total space of the vector bundle
    $\mathcal Hom(V/H,L)=\mathcal S_m\otimes\mathcal Q_{n-m}^\vee$ over $Q$. So $f$
    is affine, and
    \[
    f_*\mathcal O_{\tilde Q}=\mathrm{Sym}^*\big(\mathcal S_m^\vee\otimes\mathcal Q_{n-m}\big)
    =\bigoplus_{j\ge0}\mathrm{Sym}^j\big(\mathcal S_m^\vee\otimes\mathcal Q_{n-m}\big).
    \]
    Because $f$ is affine, $H^i(\tilde Q,\mathcal F)=H^i(Q,f_*\mathcal F)$, and
    we can apply the projection formula. Then
    \[
    H^i\big(\tilde Q,\,f^*(\Sigma^\lambda\mathcal S_m\otimes\Sigma^\mu\mathcal Q_{n-m}^\vee)\big)
    =\bigoplus_{j\ge0}H^i\big(Q,\,\Sigma^\lambda\mathcal S_m\otimes\Sigma^\mu\mathcal Q_{n-m}^\vee
    \otimes\mathrm{Sym}^j(\mathcal S_m^\vee\otimes\mathcal Q_{n-m})\big).
    \]
    Thus it suffices to show that each summand on the right vanishes when
    $|\lambda|+|\mu|=i>0$, which follows from Lemma \ref{bwbcor}.
\end{proof}

\begin{prop}\label{coh}
    We define the vector bundle $\mcM=\mcS_m^\vee\boxtimes\mcQ_{n-m}$ on $Q\times\Gr_k\times\Gr_l$. Then for all $i\geq0$,
    \[H^{-i}(\phi^*M_k^\bullet)=\wedge^i(\mcM|_{X}\oplus\mcM|_X).\]
\end{prop}

\begin{proof}
    Recall that the complex $M^\bullet$ is the tensor product of complexes $M_k^\bullet,M_l^\bullet$. Using the same notation as in the proof of \ref{main}, there are vector bundles $E_{k,1}$ and $E_{k,2}$ on $\Gr_m\times\Gr_k$ and $\Gr_{n-m}\times\Gr_k$, respectively. We denote the Koszul complex of $E_{k,1}$ (described in (\ref{koszul1})) by $M_{k,1}^\bullet$ and the Koszul complex of $E_{k,2}$ by $M_{k,2}^\bullet$. $M_k^\bullet$ is the tensor product of $M_{k,1}^\bullet$ and $M_{k,2}^\bullet$.

    We observe that it suffices to show that \[
    H^{-i}(M_k^\bullet)=\wedge^i(\mcM|_{X_k}).
    \] Supposing that this is true, we know the same for $M_l^\bullet$, and using the Kunneth formula, we obtain the cohomology of $M^\bullet$.
    
    We now compute $H^{-i}(\phi^*M_k^\bullet)$. Let $\iota_1,\iota_2$ be the inclusions of $Y_1$ and $Y_2$ (defined in (\ref{ydef})) into $Q\times\Gr_k$. We observe that since $M_{k,1}^\bullet$ and $M_{k,2}^\bullet$ are resolutions for $\mcO_{S_{m,k}}$ and $\mcO_{S_{n-m,k}}$, we can identify their derived tensor product
    \[
    M_{k,1}^\bullet\otimes^\mathbf{L}M_{k,2}^\bullet=M_{k,1}^\bullet \otimes^{\mathbf L}(\iota_2)_*(\mcO_{Y_2}).
    \] Then since tensoring a vector bundle by $(\iota_2)_*\mcO_{Y_2}$ is the same as restricting to $Y_2$ and then taking the pushforward along $\iota_2$, we have \[
    M_{k,1}^\bullet\otimes^\mathbf{L}M_{k,2}^\bullet=(\iota_2)_*(M_{k,1}^\bullet|_{Y_2}).
    \] Now it suffices to compute the cohomology of the complex \[
    E=0\to\wedge^r(E_{k,1}|_{Y_2}^\vee)\to\ldots\to \wedge^2(E_{k,1}|_{Y_2}^\vee)\to (E_{k,1}|_{Y_2}^\vee)\to \mcO_{Y_2}.
    \] 

    On $Y_2$, define vector bundles $\mcN=\mcS_{m}\boxtimes \mcQ_{n-m}^\vee$, which has fibers isomorphic to $\Hom(V/H,L)$, and $\mcF=(\mcS_{n-m}/\mcS_k)^\vee\boxtimes\mcS_{m}$, which has fibers isomorphic to $\Hom(H/V_k,L)$. There is a natural inclusion $\mcN\to E':=E_{k,1}|_{Y_2}^\vee$ which gives a short exact sequence \[
    0\to\mcN\to E'\to\mcF\to0.
    \] Such a short exact sequence gives rise to exact sequences \[
    \to S^j\mcN\to S^{j-1}\mcN\otimes E'\to S^{j-2}\mcN\otimes\wedge^2 E'\to \ldots\to \wedge^jE'\to\wedge^j\mcF\to0
    \] for all $j\geq0$. These yield a double complex $D_0$ of the form 
\[\begin{tikzcd}[cramped]
	\ldots& {S^3\mcN} &&& \\
	\ldots& {S^2\mcN\otimes E'} & {S^2\mcN} \\
	\ldots & {\mcN\otimes\wedge^2E'} & {\mcN\otimes E'} & \mcN \\
	\ldots& {\wedge^3E'} & {\wedge^2E'} & {E'} & {\mcO_{Y_2}} \\
	\ldots& {\wedge^3\mcF} & {\wedge^2\mcF} & \mcF & {\mcO_{Y_2}}
	\arrow[from=1-2, to=2-2]
	\arrow[from=2-2, to=2-3]
	\arrow[from=2-2, to=3-2]
	\arrow[from=2-3, to=3-3]
	\arrow[from=3-2, to=3-3]
	\arrow[from=3-2, to=4-2]
	\arrow[from=3-3, to=3-4]
	\arrow[from=3-3, to=4-3]
	\arrow[from=3-4, to=4-4]
	\arrow[from=4-2, to=4-3]
	\arrow[from=4-2, to=5-2]
	\arrow[from=4-3, to=4-4]
	\arrow[from=4-3, to=5-3]
	\arrow[from=4-4, to=4-5]
	\arrow[from=4-4, to=5-4]
	\arrow[from=4-5, to=5-5]
	\arrow[from=5-2, to=5-3]
	\arrow[from=5-3, to=5-4]
	\arrow[from=5-4, to=5-5]
\end{tikzcd}\]
with $H^\bullet(\Tot D_0)=0$. We observe that the bottom row is a Koszul complex for the vector bundle $\mcF^\vee$: the rank of $\mcF$ is $m(n-m-k)$, which is also the codimension of $X_k$ inside $Y_2$. There is a global section of $\mcF^\vee$ given by $$L\hookrightarrow H\to H/V_k,$$ which is zero exactly when $L$ lies inside $V_k$. Thus the bottom row has cohomology only in degree 0, where the cohomology is $\mcO_{X_k}$. Thus we can work with a double complex $D$ which removes the bottom row from $D_0$, and $H^\bullet(\Tot D)=\mcO_{X_k}$. The rows of $D$ are the complexes $S^j\mcN\otimes E[j]$ for $j\geq0$.

The spectral sequence associated to the double complex $D$ is such that \[
E_1^{p,q}=H^{q}(S^{-p}\mcN\otimes E[-p]).
\] It converges to $H^{p+q}(\Tot D)$, so for any $p,q$ where $p+q\neq0$, we have $E^{p,q}_\infty=0$. We note that $E_1^{p,q}$ is 0 for any $p>0$ or $p<q$.
Taking $q=-1$, the $E_1$ page has maps \[
0=E_1^{-2,-1}\to E_1^{-1,-1}\to E_1^{0,-1}\to E_1^{1,-1}=0.
\] Then we see from the $E_2$ page that $E^{-1,-1}_2=E_\infty^{-1,-1}=0$ and $E^{0,-1}_2=E_\infty^{0,-1}=0$, so we must have $E^{-1,-1}_1=E_1^{0,-1}$. We know that \[
E_1^{-1,-1}=H^{-1}(\mcN\otimes E[1])=\mcN\otimes H^0(E)=\mcN\otimes\mcO_{X_k}.
\] We induct on $q$ to show that $E_1^{0,q}=\wedge^q\mcN\otimes\mcO_{X_k}$. On the $E_1$ page, we have maps \[
0\to H^q(S^{-q}\mcN\otimes E[-q])\to\ldots\to H^q(\mcN\otimes E[1])\to H^q(E)\to0.
\] The inductive assumption implies that this is \[
0\to S^{-q}\mcN\otimes\mcO_{X_k}\to S^{-q-1}\mcN\otimes\wedge^1\mcN\otimes\mcO_{X_k}\to\ldots \to\mcN\otimes\wedge^{q-1}\mcN\otimes\mcO_{X_k}\to H^q(E)\to 0.
\] 
This implies that $$H^q(E)=\wedge^{-q}\mcN\otimes\mcO_{X_k},$$ which is $\wedge^{q}\mcM|_{X_k}$, as desired.
\end{proof}

\begin{prop}\label{vanishing2}
    For all $i>0,$
    \begin{align*}
        R^{i+1}\varpi_*(H^{-i}(\phi^*M^\bullet))&=0.
    \end{align*}
\end{prop}

\begin{proof}
    Since $\overline O_m$ is affine, we may further push forward from $\overline O_m\times\Gr_k\times\Gr_l$ to $\Gr_k\times\Gr_l$ and work in a fiber over some point: 
    the relevant sheaves are vector bundles over $X$, so they are flat over $X$ and therefore over $\Gr_k\times\Gr_l$. Using Grothendieck's theorem (\cite[III$_1$, Proposition~3.2.1]{EGAIII1}) for proper base change, it suffices to check our condition on each fiber.
     Fix a point $(U,W)\in\Gr_k\times\Gr_l$. Now the condition $L\subset U,W\subset H$ on $Q$ is equivalent to picking $L,H$ so that $L\subset U\cap W$ and $U+W\subset H$. Thus we can define the product of two Grassmannians \begin{align*}
        P=\Gr(L,U\cap W)\times\Gr(H/(U+W),V/(U+W)).
    \end{align*}
    and consider a commutative diagram \[
    \begin{tikzcd}
        &P\arrow[hookrightarrow,r]\arrow[d]&X\arrow[d,"h"]\\
        &Q\times\{U\}\times \{W\}\arrow[d]\arrow[hookrightarrow,r]&Q\times\Gr_k\times\Gr_l
        \arrow[d]\\
        &\{U\}\times\{W\}\arrow[hookrightarrow,r]&
        \Gr_k\times\Gr_l
    \end{tikzcd}
    \] 
    Let $\mcS_1$ be the tautological bundle on $\Gr(L,U\cap W)$ and $\mcQ_2$ the quotient bundle on $$\Gr(H/(U+W),V/(U+W)).$$ Thus it suffices to show that \begin{align*}
        H^{i+1}(P, \wedge^i((\mcS_1\boxtimes \mcQ_2^\vee)^{\oplus 2})\otimes \Sym^*(\mcS_1^\vee\boxtimes \mcQ_2))=0.
    \end{align*}

    We have \begin{align*}
        \mcS_1\boxtimes \mcQ_2^\vee\oplus\mcS_1\boxtimes \mcQ_2^\vee= (\mcS_1\boxtimes \mcQ_2^\vee)\otimes \underline{K}^{\oplus 2},
    \end{align*}     
    where $\underline{K}^{\oplus 2}$ is the trivial vector bundle of rank 2. Using the Schur functor decomposition of the exterior power, we have \begin{align*}
        \wedge^i(\mcS_1\boxtimes \mcQ_2^\vee\oplus\mcS_1\boxtimes \mcQ_2^\vee)= \bigoplus\Sigma^\lambda(\mcS_1\boxtimes \mcQ_2^\vee)\otimes\Sigma^{\lambda^*}(\underline{K}^{\oplus 2}),
    \end{align*} where $|\lambda|=i$ and $\lambda$ has rows of length at most 2 (otherwise $\Sigma^{\lambda^*}\underline{K}^{\oplus 2}=0$).
    Thus it suffices to prove that for such $\lambda$ we have 
    \begin{align*}
        H^{i+1}(P, \Sigma^{\lambda}((\mcS_1\boxtimes \mcQ_2^\vee)^{\oplus 2})\otimes \Sym^j(\mcS_1^\vee\boxtimes \mcQ_2))=0.
    \end{align*} 
    By Pieri's formula, we have 
    \begin{equation}
        \Sigma^{\lambda}(\mcS_1\boxtimes \mcQ_2^\vee)\otimes \Sym^j(\mcS_1^\vee\boxtimes \mcQ_2) \cong \oplus_{\lambda^{\prime}} \Sigma^{\lambda^{\prime}}(\mcS_1\boxtimes \mcQ_2^\vee),
    \end{equation}
    where $\lambda^{\prime} \in \mathbb{Z}^{\rk(\mcS_1\boxtimes \mcQ_2^\vee)}$ satisfy 
    \begin{enumerate}
        \item $\lambda_s \geq \lambda^{\prime}_s \geq \lambda_{s+1}$ for $1 \leq s < \rk(\mcS_1\boxtimes \mcQ_2^\vee)$;
        \item $\lambda^{\prime}_{\rk(\mcS_1\boxtimes \mcQ_2^\vee)} \leq \lambda_{\rk(\mcS_1\boxtimes \mcQ_2^\vee)}$;
        \item $\Sigma_l \lambda^{\prime}_l = i - j$.
    \end{enumerate}

    So we need to prove $$H^{i+1}(P, \Sigma^{\lambda^{\prime}}(\mcS_1\boxtimes \mcQ_2^\vee))=0.$$
    However, this follows by Theorem \ref{bwb} and Proposition \ref{kroneckergenweights}.

\end{proof}

\section{Kronecker Coefficients}\label{sec4}
In this section, we derive some combinatorial results which will be useful for studying the cohomology of certain line bundles.
Simple representations of the symmetric group $S_i$ are in bijective correspondence with partitions $\mu$ of $i$. We denote these as $V_\mu$.
We may decompose a tensor product \[
V_{\mu_1}\otimes V_{\mu_2}=\bigoplus_{\lambda\vdash i}g(\mu_1,\mu_2,\lambda)V_\lambda
\] as a sum of simple representations. The coefficients $g(\mu_1,\mu_2,\lambda)$ are called Kronecker coefficients.

\subsection{Generalities on Kronecker coefficients.} We first derive some general facts about Kronecker coefficients.

\begin{lem}\label{transpose}
    For partitions $\mu_1,\mu_2,\lambda\vdash i$, \[
    g(\mu_1,\mu_1,\lambda)=g(\mu_1,\mu_2^*,\lambda^*),
    \]
    where $\lambda^*$ is the conjugate partition.
\end{lem}
\begin{proof}
    Let $\chi^\mu$ be the character of $V_\mu$. A formula for Kronecker coefficients is \[
    g(\mu_1,\mu_2,\lambda)=\frac{1}{i!}\sum_{\pi\in S_i}\chi^{\mu_1}(\pi)\chi^{\mu_2}(\pi)\chi^{\lambda}(\pi).
    \] \cite[Chapter~I.7]{mac}. For any partition $\mu$, $\chi^{\mu^*}(\pi)=(-1)^{l(\pi)}\chi^\mu(\pi)$. Thus we may replace $\chi^{\mu_2}(\pi)\chi^{\lambda}(\pi)$ with $\chi^{\mu_2^*}(\pi)\chi^{\lambda^*}(\pi)$.
\end{proof}

We denote by $c(\alpha,\beta,\gamma)$ the Littlewood-Richardson coefficients, defined by the decomposition of Schur functors \[
\Sigma^\alpha V\otimes\Sigma^\beta V=\bigoplus_{\gamma}c(\alpha,\beta,\gamma)\Sigma^\gamma V.
\] We may relate these to Kronecker coefficients by the following.
\begin{lem}\label{LR}
    Let $j=\lfloor\frac{i}{2}\rfloor$ and $\mu_1,\mu_2\vdash i$. Then \[
    \sum_{a=0}^jg(\mu_1,\mu_2,(i-a,a))=\sum_{\substack{\alpha\vdash i-j\\\beta\vdash j}}c(\alpha,\beta,\mu_1)c(\alpha,\beta,\mu_2).
    \]
\end{lem}
\begin{proof}
    A representation $V_{\mu}$ of $S_i$ may be restricted to the subgroup $S_{i-j}\times S_j$. This restriction can be decomposed as \[
    \Res^{S_i}_{S_{i-j}\times S_j}V_\mu=\bigoplus_{\substack{\alpha\vdash i-j\\\beta\vdash j}}c(\alpha,\beta,\mu)(V_\alpha\boxtimes V_\beta).
    \] (See \cite[Chapter~4]{fh})
    We use this to compute the inner product \[
    \left<\Res^{S_i}_{S_{i-j}\times S_j}\chi^{\mu_1},\Res^{S_i}_{S_{i-j}\times S_j}\chi^{\mu_2}\right>=\sum_{\substack{\alpha\vdash i-j\\\beta\vdash j}}c(\alpha,\beta,\mu_1)c(\alpha,\beta,\mu_2).
    \] We also compute the inner product using the isomorphism between class functions of $S_i$ and the degree $i$ component of the ring of symmetric functions $\Lambda$. We refer to\cite[Chapter~I.7]{mac} for details of this isomorphism. Under this correspondence, the Schur polynomial $s_\mu$ represents the character $\chi^\mu$ and the complete homogeneous symmetric function $h_i$ represents the trivial character. Furthermore, inducing a representation of $S_{i-j}\times S_j$ to $S_i$ corresponds to multiplying the symmetric functions.
    Denote by $*$ the internal product in $\Lambda$, which corresponds to the tensor product of characters of $S_i$. Then we have the Hall inner product \begin{align*}
    \left<s_{\mu_1}*s_{\mu_2},h_{i-j}h_j\right>_\Lambda=\left<\chi^{\mu_1}\otimes\chi^{\mu_2},\Ind^{S_i}_{S_{i-j}\times S_j}(\mathbf 1_{i-j}\boxtimes \mathbf 1_j)\right>_{S_i}
    \end{align*} using the correspondence. Then by Frobenius reciprocity this is \[
    \left<\Res^{S_i}_{S_{i-j}\times S_j}(\chi^{\mu_1}\otimes\chi^{\mu_2}), \mathbf 1_{i-j}\boxtimes \mathbf 1_j\right>_{S_{i-j}\times S_j}=\left<\Res^{S_i}_{S_{i-j}\times S_j}\chi^{\mu_1}, \Res^{S_i}_{S_{i-j}\times S_j}\chi^{\mu_2}\right>_{S_{i-j}\times S_j}.
    \] Now we have \[
    \left<s_{\mu_1}*s_{\mu_2},h_{i-j}h_j\right>_\Lambda=\sum_{\substack{\alpha\vdash i-j\\\beta\vdash j}}c(\alpha,\beta,\mu_1)c(\alpha,\beta,\mu_2).
    \] Since $h_i$ is the Schur polynomial $s_{(i)}$, Pieri's formula implies that \[
    h_{i-j}h_j=\sum_{a=0}^js_{(i-a,a)}.
    \]
    Then since $s_{\mu_1}*s_{\mu_2}$ decomposes in terms of the Kronecker coefficients as with the tensor product, the multiplicity of $s_{(i-a,a)}$ is $g(\mu_1,\mu_2,(i-a,a))$. Thus we have the desired equality.
\end{proof}

We now derive necessary conditions for a Littlewood-Richardson coefficient to be nonzero. 
\begin{lem}\label{kyfan}
    If $c(\alpha,\beta,\gamma)>0$, then for all $j\geq1$, \[
    \gamma_1+\ldots+\gamma_j\leq (\alpha_1+\ldots+\alpha_j)+(\beta_1+\ldots+\beta_j)\leq \gamma_1+\ldots\gamma_{2j}.
    \]
\end{lem}
\begin{proof}

    By \cite[Theorem~11]{ful}, if $c(\alpha,\beta,\gamma)>0$, there exist Hermitian matrices $A,B$ whose eigenvalues in decreasing order are $\alpha$ and $\beta$, respectively, such that the eigenvalues of $A+B$ are $\gamma$. The lower bound then follows from \cite[Theorem~III.4.1]{bha}.
    Fixing an $n$-dimensional vector space $V$, for any Hermitian matrix $X$ with eigenvalues $\mu$ in decreasing order and any $j$, by the Ky Fan's maximum principle (see \cite[Chapter~II.1]{bha}), we have \[
    \mu_1+\ldots+\mu_j=\max_{\substack{U\subset V\\\dim U=j}}\Tr(X|_U),
    \] with $U$ running over vector subspaces of $V$ of dimension $j$. Let $U,W$ be subspaces attaining this maximum for the matrices $A$ and $B$ as above. Then for all $j$, \begin{align*}
        (\alpha_1+\ldots+\alpha_j)+(\beta_1+\ldots+\beta_j)&=\Tr(A|_U)+\Tr(B|_W)\\
        &\leq\Tr(A|_{U+W})+\Tr(B|_{U+W})\\
        &=\Tr((A+B)|_{U+W}).
    \end{align*} Since $U+W$ has dimension at most $2j$, $\Tr((A+B)|_{U+W})$ is at most the upper bound.
\end{proof}

\subsection{The vanishing result.} Now we proceed to derive necessary conditions for nonzero Kronecker coefficients.
\begin{prop}\label{kronecker}
    Let $\mu_1,\mu_2, \lambda$ be partitions of $i$ where $g(\mu_1,\mu_2,\lambda)>0$. Let $\rho=(n-1,n-2,\ldots,0)$ be the half sum of all positive roots of $\GL_n$. For a Young diagram $\omega$ with $r$ rows, we denote by $s(\omega)$ the integer sequence of length $n$ \[
    (-\omega_r,\ldots,-\omega_1,0,\ldots,0)+\rho.
    \]
    If there exist permutations $\sigma_1$ and $\sigma_2$ putting $s(\mu_1)$ and $s(\mu_2)$ into strictly decreasing order and $l(\sigma_1)+l(\sigma_2)=i+1$, then $\lambda_1>2$.
\end{prop}
\begin{proof}
    We view $\mu_1,\mu_2$ as partitions of $i$ with $r$ rows, some rows possibly being empty. We first obtain some results for $s(\mu)$ for a fixed partition $\mu=(\mu_1,\ldots,\mu_r)$ of $i$.
    Let $d=n-r$. The sequence $s(\mu)$ is some permutation of the values $1,2,\ldots,d,f(1),\ldots,f(r)$, where $f(m)=d+m-\mu_m$ is a strictly increasing function in $m$. We may eliminate from consideration any $\mu$ where $s(\mu)$ has a repeated value, since no permutation can put it in strictly decreasing order. This occurs exactly when one of the values $f(m)$ conflicts with a value in $1,\ldots,d$, i.e. \[
    1\leq d+m-\mu_m\leq d
    \] for some $m$. We proceed assuming that no such $m$ exists. 

    Thus we can divide the integers $1,\ldots,r$ into two groups: either $f(m)\leq 0$, or $f(m)>d$. We denote by $L=L(\mu)$ the number of indices such that $f(m)\leq0$, i.e. $\mu_m\geq d+m$. Call a pair of indices $i<j$ an \textit{inversion} in $s(\mu)$ if $s(\mu)_i<s(\mu)_j$ (so that these values are not in decreasing order). The number of inversions is equal to the length of the permutation $\sigma$ putting $\mu$ in decreasing order. We see that $s(\mu)$ is organized as follows: \[
    \underbrace{f(r),f(r-1),\ldots,f(L+1)}_{> d},\underbrace{f(L),\ldots,f(1)}_{\leq0},d,d-1,\ldots,1.
    \] Inversions occur between the middle and right blocks, so the number of inversions is $Ld$.

    By definition, $\mu_L\geq d+L$. Since $\mu$ is decreasing, \begin{align}\label{bound}
    \mu_1+\ldots+\mu_L\geq L(d+L)=Ld+L^2.
    \end{align}

    Now suppose a partition $\lambda$ is such that $\lambda_1\leq 2$. We wish to show that $g(\mu_1,\mu_2,\lambda)=0$ for any $\mu_1,\mu_2$ with $l(\sigma_1)+l(\sigma_2)=i+1$. By Lemma \ref{transpose}, this is equivalent to $g(\mu_1,\mu_2^*,(i-a,a))=0$ for all two-row partitions $(i-a,a)$ of $i$.
    By Lemma \ref{LR} we have \[
    \sum_{a=0}^jg(\mu_1,\mu_2^*,(i-a,a))=\sum_{\substack{\alpha\vdash i-j\\\beta\vdash j}}c(\alpha,\beta,\mu_1)c(\alpha,\beta,\mu_2^*).
    \] Since all summands on the left are nonnegative, it suffices to show that the right side is 0. 
    
    Suppose that there exist $\alpha,\beta$ such that $c(\alpha,\beta,\mu_1),c(\alpha,\beta,\mu_2^*)>0$. The lower bound of Lemma \ref{kyfan} for $\mu_1$ with $j=L(\mu_1)=L$ yields \[
    \mu_{1,1}+\ldots+\mu_{1,L}\leq \alpha_1+\ldots+\alpha_L+\beta_1+\ldots+\beta_L.
    \] Equation \ref{bound} gives \[
    \mu_{1,1}+\ldots+\mu_{1,L}\geq l(\sigma_1)+L^2=i+1-l(\sigma_2)+L^2.
    \] Thus \begin{align}\label{ineq1}
    \alpha_1+\ldots+\alpha_L+\beta_1+\ldots+\beta_L\geq i-l(\sigma_2)+L^2+1.
    \end{align}

    Now the upper bound of Lemma \ref{kyfan} for $\gamma=\mu_2^*$ yields \[
    \alpha_1+\ldots+\alpha_L+\beta_1+\ldots+\beta_L\leq(\mu_2^*)_1+\ldots+(\mu_2^*)_{2L}.
    \] Viewing $\mu_2$ as a Young diagram, the quantity on the right is the number of boxes in the first $2L$ columns. Let $L'=L(\mu_2)$. By definition, $$\mu_{2,L'}\geq d+L',$$
    where $d$ is the difference between $n$ and the length of $\mu_2$. Thus the first $L'$ rows of $\mu_2$ each contain at least $\max(d+L'-2L,0)$ boxes not in the first $2L$ columns, so \[
    i-(\mu_2^*)_1+\ldots+(\mu_2^*)_{2L}\geq L'\cdot\max(d+L'-2L,0).
    \] Since $d+L^\prime \geq 2L$ we get  \begin{align}
        (\mu_2^*)_1+\ldots+(\mu_2^*)_{2L}\leq i-L'd+L^2.
    \end{align} But $L'd=l(\sigma_2)$, which contradicts Equation \ref{ineq1}. Thus the Kronecker coefficients $g(\mu_1,\mu_2,\lambda)=0$ vanish.
\end{proof}

\subsection{The case of dominant integral weights.} In this subsection we generalize Proposition \ref{kronecker} from the case of partitions to the case of non-increasing sequences of (possibly negative) integers.

\begin{notn}
    For positive $a, b \in \mathbb{Z}$ and $$\mu_1 \in \mathbb{Z}^a,\mu_2 \in \mathbb{Z}^b, \lambda \in \mathbb{Z}^{ab},$$
    define $g^{a, b}(\mu_1,\mu_2,\lambda)$ as the multiplicity of the summand $$\Sigma^{\mu_1}(F_a) \otimes \Sigma^{\mu_2}(F_b)$$ in $\Sigma^{\lambda}(F_a \otimes F_b)$, where $F_a$ and $F_b$ are $K$-vector spaces of dimension $a$ and $b$, respectively.
\end{notn}
By Schur–Weyl duality, when $\mu_i$ and $\lambda$ are partitions, $g^{a, b}(\mu_1,\mu_2,\lambda)$ coincides with the standard notion of Kronecker coefficients.

\begin{notn}
    For $\lambda \in \mathbb{Z}^{l}$, set $$\lambda_+=(\max(\lambda_k, 0))_k;$$ 
$$\lambda_{-} = -(\min(\lambda_k, 0))_k.$$
\end{notn}

\begin{prop}\label{kroneckergenweights}
    Fix positive $a, b \in \mathbb{Z}$. Let $$\mu_1 \in \mathbb{Z}^a,\mu_2 \in \mathbb{Z}^b, \lambda \in \mathbb{Z}^{ab}$$ be non-increasing sequences. Assume $g^{a, b}(\mu_1,\mu_2,\lambda)>0$. 
    If there exist permutations $\sigma_1$ and $\sigma_2$ putting $s(\mu_1)$ and $s(\mu_2)$ into strictly decreasing order and $l(\sigma_1)+l(\sigma_2) \geq |\lambda|_+ + 1$, then $\lambda_1>2$.
\end{prop}

We will reduce the statement of Proposition \ref{kroneckergenweights} to Proposition \ref{kronecker}. The following lemma will be useful for this.

\begin{lem}\label{reduction}
    In the setting of Proposition \ref{kroneckergenweights}, there exist $\alpha_1 \in \mathbb{Z}^a$, $\alpha_2 \in \mathbb{Z}^b$ partitions of $|\lambda_+|$ such that 
    \begin{enumerate}
        \item $\alpha_{1, k} \geq \mu_{1, k}$ for all $k$,
        \item $\alpha_{2, l} \geq \mu_{2, l}$ for all $l$,
        \item $g(\alpha_1, \alpha_2, \lambda_+) >0$.
    \end{enumerate}
\end{lem}
\begin{proof}
    Indeed, we have 
    \begin{equation}\label{1}
        \Sigma^{\lambda}(F_a \otimes F_b) \hookrightarrow \Sigma^{\lambda_+}(F_a \otimes F_b) \otimes \Sigma^{\lambda_-}(F_a \otimes F_b)^\vee.
    \end{equation}
    Decomposing 
    \begin{equation}
        \Sigma^{\lambda_+}(F_a \otimes F_b) \cong \oplus (\Sigma^{\alpha_1}(F_a) \otimes \Sigma^{\alpha_2}(F_b))^{\oplus g(\alpha_1, \alpha_2, \lambda_+)};
    \end{equation}
    \begin{equation}
        \Sigma^{\lambda_-}(F_a \otimes F_b) \cong \oplus (\Sigma^{\beta_1}(F_a) \otimes \Sigma^{\beta_2}(F_b))^{\oplus g(\beta_1, \beta_2, \lambda_-)},
    \end{equation}
    we get that $$\Sigma^{\lambda_+}(F_a \otimes F_b) \otimes \Sigma^{\lambda_-}(F_a \otimes F_b)^\vee$$ decomposes into summands of the form 
    $$(\Sigma^{\alpha_1}(F_a) \otimes \Sigma^{\beta_1}(F_a)^\vee) \otimes (\Sigma^{\alpha_2}(F_b) \otimes \Sigma^{\beta_2}(F_b)^\vee)$$
    with multiplicities $g(\alpha_1, \alpha_2, \lambda_+)g(\beta_1, \beta_2, \lambda_-)$. 

    By \ref{1}, at least one of the multiples is non-zero. Pick these $\alpha_1$, $\alpha_2$, $\beta_1$, $\beta_2$. We know that 
    \begin{equation}
        \Sigma^{\mu_1}(F_a) \hookrightarrow (\Sigma^{\alpha_1}(F_a) \otimes \Sigma^{\beta_1}(F_a)^\vee);
    \end{equation}
    \begin{equation}
        \Sigma^{\mu_2}(F_b) \hookrightarrow (\Sigma^{\alpha_2}(F_b) \otimes \Sigma^{\beta_2}(F_b)^\vee).
    \end{equation}
    This implies that for large enough integer $N$,
    we have $$c(\mu_1  + N, \beta_1,\alpha_1 + N) >0,$$
    and similarly for the triple $\mu_2$, $\beta_2$, $\alpha_2$. The assertion follows.
\end{proof}

\begin{proof}[Proof of Proposition \ref{kroneckergenweights}]
    By Lemma \ref{LR}, since $\lambda_+$ has at most 2 columns, we get that there exist partitions $\nu, \tau$ of $|\lambda_+|$ such that $$c(\nu, \tau, \alpha_1), c(\nu, \tau, \alpha_2^*)>0.$$
    By Lemma \ref{kyfan} we have
    \begin{equation}
        \alpha_{1,1}+\ldots+\alpha_{1,L}\leq (\nu_1+\ldots+\nu_L)+(\tau_1+\ldots+\tau_L)\leq (\alpha^*_{2})_{1}+\ldots+(\alpha^*_{2})_{2L}.
    \end{equation}
    Keep the notation from Proposition \ref{kronecker}: let $d=n-b$, $L=L(\mu_1)$ and $L^{\prime} = L(\mu_2)$. (Note that the notion $L(\alpha)$ makes sense for general non-increasing sequence of integers, and the estimates in Proposition \ref{kronecker} also hold for such $L(\alpha)$). 
    
    By property (1) in Lemma \ref{reduction} we get 
    \begin{equation}
        \alpha_{1,1}+\ldots+\alpha_{1,L}\geq \mu_{1,1}+\ldots+\mu_{1,L}\geq l(\sigma_1)+L^2.
    \end{equation}
    Moreover, as in the proof of Proposition \ref{kronecker}, we get that 
    \begin{equation}
        (\alpha^*_{2})_{1}+\ldots+(\alpha^*_{2})_{2L} \leq |\lambda_+| - L'\cdot\max(d+L'-2L,0) \leq |\lambda_+| - L'(d+L'-2L).
    \end{equation}
    Hence $$l(\sigma_1)+l(\sigma_2) \leq |\lambda_+| - (L-L^\prime)^2,$$
    hence the result.
\end{proof}
\newpage

\printbibliography

\end{document}